\documentclass[12pt]{article}
\usepackage{amssymb,amsmath}

\def\a{\alpha}
\def\b{\beta}
\def\s{\sigma}
\def\frak{\mathfrak}
\def\Cal{\mathcal}

\begin{document}

\makeatletter	   
\renewcommand{\ps@plain}{%
     \renewcommand{\@oddhead}{\textrm{\scriptsize FAR EAST JOURNAL OF DYNAMICAL SYSTEMS,
 vol. 4 (2), pp. 125-134, 2003.}\hfil\textrm{}}%
     \renewcommand{\@evenhead}{\@oddhead}%
     \renewcommand{\@oddfoot}{}
     \renewcommand{\@evenfoot}{\@oddfoot}}
\makeatother     

\title {\sl Chaos in Topological Spaces\\
\vskip 10pt
}
 \author{John Taylor\\
 Department of Mathematics\\
University of Colorado, Boulder\\
Boulder, CO USA\\
{\small jtaylor@euclid.colorado.edu}
\thanks{I thank Paul Meyer for his helpful suggestions.  I also thank Erwin Kronheimer and Mati Rubin for the very useful conversations.}
}
\date{} 
\maketitle

\begin{abstract}

\noindent We give a definition of chaos for a continuous self map $f$ of a general
topological space.  This definition coincides with Devanney's definition
for chaos when the topological space happens to be a metric space.  This generaization is nescessry
to the study chaotic phenomena in relation to topological entropy.

\noindent We show that in a uniform Hausdorff space, there is
a meaningful definition of sensitive dependence on initial
conditions, and prove that if $f$ is chaotic on a such a
space, then $f$ necessarily has sensitive dependence on initial conditions.
The proof is interesting in that it explains very clearly what causes a chaotic
process to have sensitive dependence.  

Finally, we construct a chaotic map on a non-metrizable topological space.

\end{abstract}

\centerline {INTRODUCTION}
\bigskip

This note is put forward as a meaningful extention of the notion of chaos from metric to topological spaces.  This generalization is worthwhile from different points of view.  For example, it provides for the direct comparison and study of possible relations between chaos and topological entropy.  Furthermore, as we show by construction, chaotic phenomena are present in non-metrizable spaces.

The point of departure here is the work of Banks et al.. \cite{1}.  This result states that, in a metric space, the first two
components of Devaney's definition \cite{2} --- both purely topological --- imply the
third:  sensitive dependence on initial conditions, which has a metric
formulation.  Because sensitive dependence is the only aspect of chaotic
behavior that is open to experimental detection, it
is considered by the scientific community to be the hallmark of chaotic
behavior.  Hence, anything one can say about sensitive dependence is important.  We constructively explore the most general setting in which this phenomena occurs (a uniform Hausdorff space) and ellucidate its cause. 
\bigskip

\centerline {SOME PRELIMINARIES}
\bigskip

Since we will be concerned with uniform topological spaces, we recall
briefly some definitions and facts concerning them.  For more details we
refer the reader to \cite{3}.

Let $X$ be a nonvoid set and $U,V$ nonvoid subsets of $X\times X$.  Then
$U$ and $V$ are relations on $X$.  If $U\subset X\times X$, then
$U^{-1} = \{ (y,x) \mid (x,y) \in U\}$; $U\circ V = \{ (x,z)\mid$ for
some $y,(x,y)\in V$ and $(y,z)\in U\}$; and the diagonal $\Delta (X) =
\Delta = \{ (x,x)\mid x\in X\}$.
\bigskip

{\bf Definition 1}  A {\bf uniformity} for a set $X$ is a nonvoid
collection $\frak U$ of subsets of $X\times X$ such that

\ (1)\ each member of $\frak U$ contains $\Delta$;
\smallskip

\ (2)\ if $U\in {\frak U}$, then $U^{-1} \in {\frak U}$;
\smallskip

\ (3)\ if $U\in {\frak U}$, then $V\circ V\subset U$ for some $V \in
{\frak U}$;
\smallskip

\ (4)\ if $U$ and $V$ are members of $\frak U$, then $U\cap V \in 
{\frak U}$; and
\smallskip

\ (5)\ if $U \in {\frak U}$ and $U\subset V\subset X\times X$, then $V \in
{\frak U}$.
\smallskip

The pair $(X,{\frak U})$ is called a {\bf uniform space}.

A subcollection $\Cal B$ of a uniformity $\frak U$ is a {\bf base} for
$\frak U$ if and only if each element of $\frak U$ contains an element of
$\Cal B$.

For $x\in X$ and $U\in {\frak U}$, define
$$U[x] = \{ y\in X \mid (x,y)\in U\} .$$
Extend this to subsets $A$ of $X$ as follows:
$$U[A] = \bigcup_{x\in A} U[x] = \{ y\in X \mid (x,y)\in U\}.$$
\bigskip

{\bf Definition 2}\  If $(X, {\frak U})$ is a uniform space, the
topology $\Cal T$ of the uniformity, called the {\bf uniform topology},
is the colleciton of all subsets $T$ of $X$ such that for each $x\in T$
there is a $U\in {\frak U}$ such that $U[x]\subset T$.

In fact, for each $x \in X$, the collection ${\Cal U}_x = \{ U[x] \mid U
\in {\frak U}\}$ forms a neighborhood base at $x$.

In this paper we employ the idea, due to A. Weil \cite{4}, of a uniform
neighborhood system.
\bigskip

{\bf Definition 3}\  A {\bf uniform neighborhood system} is a nonvoid
set $X$, an arbitrary index set $\Lambda$, an ordering $\succcurlyeq$ in $\Lambda$,
and a function $V : \Lambda
\times X\longrightarrow \wp (X)$ such that the following conditions are
satisfied:

\ (1)\ for all $\alpha \in \Lambda$ and all $x \in X$, $V_\alpha (x)$ is a
subset of $X$ to which $x$ belongs;
\smallskip

\ (2)\ the relation $\succcurlyeq$ directs the set $\Lambda$;
\smallskip

\ (3)\ if $\alpha \succcurlyeq \beta$, then $V_\alpha (x)\subset V_\beta
(x)$ for all $x$.
\smallskip

\ (4)\ for all $\alpha \in \Lambda$ there is a $\beta \in \Lambda$ such
that $y \in V_\alpha (x)$ whenever $x \in V_\beta (y)$; and
\smallskip

\ (5)\ for all $\alpha \in \Lambda$ there is a $\beta \in \Lambda$ such
that $z \in V_\alpha (x)$ whenever $y \in V_\beta (x)$ and $z \in V_\beta
(y)$.

If $(V,\succcurlyeq)$ is a uniform neighborhood system for $X$, the
collection of all sets $\{ (x,y) \mid y \in V_\alpha (x)\}$ is base for a
uniformity $\frak U$ of $X$, called the uniformity of
the system.  This uniformity has the property that for each $\alpha \in
\Lambda$, there is some $U \in {\frak U}$ such that $U[x] \subset
V_\alpha (x)$ for all $x \in X$.  Thus a uniform neighborhood system
forms a neighborhood base at $x$, making $X$ a topological space.
\bigskip

\centerline {MAIN SECTION}
\bigskip

For the remainder of the note, unless otherwise specified, $X$ will 
represent an infinite uniform Hausdorff space.  Also, we refer to the
five axioms of a uniform neighborhood system found in definition 3 by
UNS$(*)$, where $*$ is a positive integer between 1 and 5 inclusive.
\bigskip
\bigskip
\bigskip

{\bf Definition 4}\  If $X$ is any nonvoid set and $f : X
\longrightarrow X$ any mapping of $X$, we define the set $O^+_f (x) = \{
x,f(x),f^2 (x),\dots \}$, where $f^n (x) = f[f^{n-1} (x)]$.  $O^+_f (x)$
is called the {\bf forward orbit} of $x$ under $f$.
\bigskip

{\bf Notation:} For simplicity, if the mapping $f$ is understood and there is no danger of confussion, we will
write $O(x)$ in place of  $O^+_f (x).$  Also, we will sometimes write $f^k(x)=x_k$, so that
$O(x)=\{x,x_1,x_2,\dots\}.$

{\bf Definition 5}\ If $X$ is any topological space and $f : X
\longrightarrow X$ any continuous mapping of $X$, then $f$ is said to be
{\bf transitive} on $X$ if for any two nonvoid open subsets $U$ and $V$
of $X$, there is a positive integer $n$ such that $f^n (U) \cap V \ne \emptyset$.
\bigskip

{\bf Definition 6}\ Let $X$ be any nonvoid set and $f : X
\longrightarrow X$ any mapping.  If for some positive integer $n$, $f^n
(x) = x$ but $f^k (x) \ne x$ for all $k$, $1 \le k \le n-1$, then $x$ is
called a {\bf periodic point of $f$} with primitive period $n$.
\bigskip

{\bf Notation:} We denote the set of all periodic points in $X$ by $Per(X).$
\bigskip

{\bf Definition 7}\ Let $X$ be a topological space and $f : X
\longrightarrow X$ a continuous map.  $f$ is said to be {\bf chaotic} on
$X$ provided that any two nonvoid open subsets of $X$ share a periodic orbit
-- that is, given $U$ and $V$ nonvoid and open in $X$, there is a periodic point $p$ in $X$ such that
$$ U\cap O (p)\ne\emptyset \quad \text{and}\quad V\cap O (p)\ne\emptyset . $$
\bigskip

{\bf Proposition}\
Let $X$ be a topological space and $f:X\longrightarrow X$ be continuous.  $f$ is chaotic on $X$ if and only if $f$ is transitive and the periodic points of $f$ are dense in $X$.
\bigskip

{\bf Proof}\ If $f$ is chaotic on $X,$ then every pair of open sets share a periodic orbit; in particular, the periodic points are dense in $X$.  Given $U,V\subset X$ non-void open sets, there is a periodic point $p$ and $x,y\in O(p)$ with $x\in U$ and $y\in V$ and some positive integer $k$ such that $f^k(x)=y,$ so that $f^k(U)\cap V\ne\emptyset.$

On the other hand, the transitivity of $f$ on $X$ implies that given $U,V\subset X$ open and non-void, there is a $k$ such that for some $x\in U$, $f^k(x)\in V.$  Now define $W=f^{-k}(V)\cap U.$  Then $W$ is open and non-void with the property that $f^k(W)\subset V.$  But since the periodic points of $f$ are dense in $X$, there is a $p\in W$ such that $f^k(p)\in V.$  This shows that $U,V$ share the orbit of $p$ --- that is, $U\cap O (p)\ne\emptyset \quad \text{and}\quad V\cap O (p)\ne\emptyset.$
\bigskip

\noindent {\bf Remark.}  Easy examples show that topological transitivity and the density of the periodic points are
independent; for example, a rotation of a cylinder by $\pi/4$ radians has dense periodic points but no transitivity, while an irrational rotation of the circle has transitivity but no periodic points.
\bigskip

{\bf Definition 8}\ Let $f : X \longrightarrow X$ be continuous.  $f$
is said to have {\bf sensitive dependence on initial conditions} provided
there exists $\alpha\in \Lambda$ such that for all $x \in X$ and every
neighborhood $N(x)$ of $x$, there exists $y \in N(x)$ ($y \ne x$), and $n
\in {\Bbb N}$ such that $f^n (x) \notin V_{\alpha}$ $(f^n (y))$.
$\alpha$ is called the {\bf sensitivity index} of $f$.
\bigskip

{\bf Lemma 1}\ There exists $\eta\in\Lambda$ such that for all
$x\in X$ there is a $q\in Per(X)$ with
$$V_{\eta} (q_k)\cap\{ x\} = \emptyset\quad\text{for all}\quad q_k \in
O(q).$$
\bigskip

\noindent {\bf Proof.}\ Let $r,s$ be distinct periodic points, so that
$O(r)\cap O(s) = \emptyset$.  First, observe that
there exists $\eta\in\Lambda$ such that for all $r_k \in O(r)$ and $s_k \in O(s)$,
$$V_{\eta} (r_k) \cap V_{\eta} (s_k) = \emptyset .$$
For $O(r)$ and $O(s)$ are finite and disjoint, so that
their union is finite and for every pair $x,y$ in the union, there is a
$\beta \in \Lambda$ such that $V_{\beta} (x) \cap V_{\beta} (y) =
\emptyset$, since $X$ is a uniform Hausdorff space.  Take the supremum over all such indices and call it $\eta$.

Now let $x \in X$.  Since $r,s$ are distinct periodic points, if $x
\notin V_{\eta} (r_k)$ for all $r_k \in O(r)$, letting $r =
q$ gives the result.  If $x \in V_{\eta} (r_k)$ for some $r_k \in
O(r)$, then $x \notin V_{\eta}(s_k)$ for all $s_k \in O(s)$ and again, letting $s = q$, we are done. 
\bigskip

{\bf Lemma 2}\ There exists $\alpha\in\Lambda$ such that for all $x \in X$
there is a periodic point $q \in X$ with
$$V_\alpha (q_k) \cap V_\alpha (x) = \emptyset \quad \text{for all} \quad
q_k \in O (q).$$
\bigskip

\noindent {\bf Proof.}\  Select any $x\in X$ and let $\eta$ and $q$ be as in Lemma 1, so that for each $q_k
\in O (q),$ $x\notin V_{\eta}(q_k).$  Using UNS(5), given this $\eta$, there is a $\b\in\Lambda$
such that
$$x \in V_{\eta}(q_k) \quad \text{whenever} \quad y \in V_\b (q_k) \quad \text{and} \quad x \in V_\b (y).$$
Thus for all $y \in V_\b (q_k)$, $x \notin V_\b (y),$ since otherwise, we would have that $x\in V_{\eta}(q_k)$, contradicting what we are assuming to be true.

\noindent Now by UNS(4), given
$\beta$, there exists $\beta_0$ such that if $y \in V_{\b_0} (x)$, then $x \in V_\b
(y)$.  But we have just shown that for all $y \in V_\b (q_k)$, $x \notin V_\b (y)$.
Therefore, for all $y \in V_\b (q_k)$, $y \notin V_{\b_0} (x)$.  Hence $V_\b (q_k) \cap V_{\b_0}
(x) = \emptyset$; and since $\succcurlyeq$ directs $\Lambda$, we may choose $\alpha \in
\Lambda$ such that $\alpha \succcurlyeq \beta$ and $\alpha \succcurlyeq \beta_0$.  Thus $\alpha$
is independent of $x$ (since $\eta$ is independent of $x$ and $\b$ depends only on $\eta,$ while $\beta_0$ depends only on $\beta)$ and we see
(by UNS(3)) that
$$V_\alpha (x) \cap V_\alpha (q_k) = \emptyset \quad \text{for all} \quad
q_k \in O (q)$$

\noindent which completes the proof.
\bigskip

{\bf Theorem 1.}\quad Let $X$ be a uniform Hausdorff space and $f : X
\longrightarrow X$ a continuous mapping.  If $f$ is chaotic, then $f$ has
sensitive dependence on inital conditions.
\bigskip

\noindent {\bf Proof.}\quad Given $x, \a, q$ as in lemma 2.

Denoting $q_i = f^i(q)$,
\smallskip

let $W_0 = V_\a(q_n),$ and for $1\le i\le n,$
\smallskip

$W_i = f^{-1}(W_{i-1})\cap V_\a(q_{n-i}).$ 
\smallskip

Then $f^i(W_n)\subset W_{n-i}\subset V_\a(q_i).$  Choose any open neighborhood $N_x$ of $x$, calling
\smallskip

$U=N_x\cap V_\a(x).$
\smallskip

Since the periodic points are dense in $X$, there exists a $p\in U$ with primitive period $n.$
\smallskip

By transitivity, there exists $k\in {\Bbb N}$ such that $f^k(U)\cap W_n\not=\emptyset,$
\smallskip

which means there exists $y\in U$ such that $f^k(y)\in W_n.$
\smallskip

Write $k=an-b$ for $0\le b\le n-1.$
\smallskip

Then $f^{an}(y)=f^b(f^k(y))\in f^b(W_n)\subset W_{n-b}\subset V_\b(q_b),$
\smallskip

and since the primitive period of $p$ is $n$, $f^{an}(p)\in U.$
\smallskip

By lemma 2  $V_\a(q_b)\cap V_\a(x)=\emptyset,$ which implies that $f^{an}(y),f^{an}(p)$ are $\a$-separated.  
\smallskip

The question is, where is $f^{an}(x)$? 
\smallskip

If $f^{an}(x)\in V_\a(q_b)$ then $f^{an}(x),f^{an}(p)$ are $\a$-separated, since $V_\a(q_b)\cap V_\a(p)=\emptyset.$ 
\smallskip

If $f^{an}(x)\not\in V_\a(q_b),$ then $f^{an}(x)\not\in V_\a(f^{an}(y)),$
\smallskip

Because, $f^{an}(y)\in V_\a(q_b)\subset V_\a(q_b).$
\smallskip

Hence, $f^{an}(x)\not\in V_\a(f^{an}(y),$ so that $f^{an}(x),f^{an}(y)$ are $\a$-separated.
\smallskip

Likewise, if $f^{an}(x)\in V_\a(x)$ then $f^{an}(x),f^{an}(y)$ are $\a$-separated.
\smallskip

In a similar way, one shows that
\smallskip

If $f^{an}(x)\not\in V_\a(x),$ then $f^{an}(x)\not\in V_\a(f^{an}(p)),$
\smallskip

Because,$V_\a(x)\supset V_\a(x)\supset U,$ and $p\in U.$
\smallskip

Hence, $f^{an}(x)\not\in V_\a(x)\Rightarrow f^{an}(x)\not\in V_\a(p),$ so that $f^{an}(x),f^{an}(p)$ are $\a$-separated.
\smallskip
In any case, $f^{an}(x)$ is $\a$-separated from one of $f^{an}(y),f^{an}(p),$
\smallskip

showing that $f$ has sensitive dependence on initial conditions, with index of sensitivity $\a.$

\bigskip

\centerline {A CHAOTIC MAP ON A NON-METRIZABLE TOPOLOGICAL SPACE}
\bigskip
 
The cantor set, as a topological space, can be realized as the set of all sequences of zeros and ones with appropriate basic neighborhoods.  That is, $C=\{f:\Cal{N}\rightarrow \{ 0,1\}\}$.

We mimic this construction with $\Cal{N}$ replaced by an uncountable set $A$ and an adjustment of basic neighborhoods.
\bigskip

The Cantor Space $C$ over $A$:
\bigskip

Let $A$ be a infinite set.  Let $C=\{0,1\}^A$ be the set of functions from $A$ to $\{0,1\}.$
\bigskip

For each finite set $B\subset A$ and each map $\phi\in\{0,1\}^B$, let $N(B,\phi)=\{f\in C | f(x)=\phi(x)\ \text{for}\ x\in B\}.$
\bigskip

Let $\Cal{T}_C$ be the topology whose basis is the set $\{N(B,\phi) | B\subset A\ \text{finite}\ \text{and}\ \phi\in\{0,1\}^B\}.$
\bigskip

$(C,\Cal{T}_C)$ is the Cantor Space which, not having a countable neighborhood base, is a non-metrizable topolgical space.
\bigskip

{\bf Lemma 3}

\smallskip
Let $A=\cup_{\a\in\Cal{A}}A_{\a}$ where $A_{\a}=\{x_{\a ,i} | i\in\Cal{N}\}$ is countable and $\a\ne\beta\implies A_{\a}\cap A_{\beta}=\emptyset.$  For each finite set $S\subset\Cal{A}$, $k\in\Cal{N},$ define $K(S,k)=\{x_{\a ,i} | \a\in S,\ 1\le i\le k\}.$  Then $\{N(K(S,k),\psi) | \psi\in \{0,1\}^K \}$ is a basis for $\Cal{T}_{\Cal{C}}.$
\bigskip

{\bf Proof}
This is an easy check.
\bigskip

{\bf Theorem 2}
If $\s :C\longrightarrow C$ is defined by $(\s f)(x_{\a ,i})=f(x_{\a ,(i+1)})$ then $\s$ is chaotic.
\bigskip

{\bf Proof}
$\s$ is topologically transitive:
\bigskip

Let $U,V\subset\Cal{T}_\Cal{C}.$  By lemma 3, there exist sets $K(S_1,k_1), K(S_2,k_2)$ and $\phi_1,\phi_2$ such that $N(K(S_1,k_1),\phi_1)\subset U$ and $N(K(S_2,k_2),\phi_2)\subset V.$
\bigskip
Let $S=S_1\cup S_2,\ K=K(S,k_1+k_2)$ and $\psi:K\longrightarrow\{0,1\}.$
\bigskip
Define $\psi$ on $K(S,k_1+k_2)$ by
$$
\psi(x_{\a ,i})= \begin{cases} \phi_1(x_{\a ,i}) &\text{if $\a\in S_1\ and\ i\le k_1$}\\ \phi_2(x_{\a ,i-k_1}) &\text{if $\a\in S_2\ and\ i > k_1$}\\ 0 &\text{otherwise} .\end{cases}
$$
Since $\psi(x)=\phi_1(x)$ for $x\in K(S_1,k_1)$ it follows that $N(K(S,k_1+k_2),\psi)\subset N(K(S_1,k_1),\phi_1).$
\bigskip
If $f\in N(K(S,k_1+k_2),\psi)$ then $(\s^{k_1}f)(x)=\phi_2(x)$ for $x\in K(S_2,k_1)$ and it follows that $\s^{k_1}N(K(S,k_1+k_2),\psi)\subset N(K(S_2,k_2),\phi_2).$
\bigskip
That is, $\s^{k_1}(U)\cap V\not=\emptyset.$
\bigskip
The set of periodic points of $\s$ are dense in C:
\bigskip

Let $U\in\Cal{T}_{\Cal{C}}.$  By lemma 3, there exists $N(K(S,k),\phi)\subset U.$
\bigskip

Define $f\in N(K(S,k),\phi)$ by
$$
f(x)= \begin{cases} \phi(x_{\a ,i}) &\text{if $x\in\{x_{\alpha, j} | \alpha\in S, 1\le i\le k, \text{and}\ j=k\lambda+i, \lambda=0,1,2\dots\}$} \\  0 &\text{otherwise} .\end{cases}
$$

Then $\s^kf=f.$ 
\bigskip

{\bf Note} $C$ admits a uniformity but is not metrizable.
\vfill\eject

\end{document}